\documentclass[10pt]{amsart}
\usepackage{graphicx,amssymb}
\usepackage{a4wide}

\theoremstyle{plain}
\newtheorem{theorem}{Theorem}[section]
\newtheorem{lemma}[theorem]{Lemma}

\newtheorem{corollary}[theorem]{Corollary}

\theoremstyle{definition}
\newtheorem{definition}[theorem]{Definition}
\newtheorem{claim}[theorem]{Claim}
\newtheorem{example}[theorem]{Example}

\def\Z{\mathbb Z}
\def\N{\mathbb N}

\def\D{\mathcal D}

\def\wedgeL{\wedge_L}

\def\veeL{\vee_{\!L}}

\let\le\leqslant
\let\ge\geqslant
\def\len{\operatorname{len}}
\def\infs{\inf{\!}_ s}
\def\sups{\sup{\!}_s}
\def\lens{\len{\!}_s}

\def\Lmax{\operatorname{L^{max}}}
\def\St{{\operatorname{St}}}
\def\t{\operatorname{\it t}}
\def\INF{\t_{\inf}}
\def\SUP{\t_{\sup}}
\def\LEN{\t_{\len}}

\begin{document}

\title
[Periodically geodesic powers in Garside groups]
{Some power of an element in a Garside group is\\ conjugate
to a periodically geodesic element}

\author{Eon-Kyung Lee \and Sang-Jin Lee}

\thanks{This work was supported by the Korea Research Foundation Grant
funded by the Korean Government (MOEHRD, Basic Research Promotion Fund)
(KRF-2006-312-C00469).}

\address{Department of Applied Mathematics, Sejong University,
    Seoul, 143-747, Korea}
\email{eonkyung@sejong.ac.kr}
\address{Department of Mathematics, Konkuk University,
    Seoul, 143-701, Korea}
\email{sangjin@konkuk.ac.kr}
\date{\today}

\begin{abstract}
We show that for each element $g$ of a Garside group,
there exists a positive integer $m$ such that
$g^m$ is conjugate to a periodically geodesic element $h$,
an element with $|h^n|_\D=|n|\cdot|h|_\D$ for all integers $n$,
where $|g|_\D$ denotes the shortest word length of $g$
with respect to the set $\D$ of simple elements.
We also show that there is a finite-time
algorithm that computes, given an element of a Garside group, its
stable super summit set.

\end{abstract}

\maketitle

\section{Introduction}

For a finitely generated group $G$ and
a finite set $X$ of semigroup generators for $G$,
the \emph{translation number} with respect to $X$ of
an element $g\in G$ is defined by
$$
\t_{G,X}(g)=\lim_{n\to \infty}\frac{|g^n|{}_X}n,
$$
where $|\cdot|_X$ denotes the shortest word length in the alphabet $X$.
If there is no confusion about the group $G$,
we simply write $\t_X(g)$ instead of $\t_{G,X}(g)$.
When $A$ is a group generator, $|g|_A$ and $\t_A(g)$ indicate
$|g|_{A\cup A^{-1}}$ and $\t_{A\cup A^{-1}}(g)$, respectively.
The following list is a part of the previous works
on the discreteness properties of translation numbers in geometric
and combinatorial groups:
\cite{Gro87,GS91,Swe95,Kap97,Bes99,Con00,CMW04,Lee06,LL06b}.
An element $g$ is said to be
\emph{periodically geodesic} with respect to $X$
if $|g^n|_X=|n|\cdot|g|_X$ for all integers $n$.

It is well-known that, in a word hyperbolic group,
all translation numbers are rational with uniformly bounded denominators.
This follows from a claim of Gromov~\cite{Gro87} and
was accurately proved by Swenson~\cite{Swe95}.
In fact, Swenson proved that for every element $g$
of a word hyperbolic group, there exists a positive integer $m$
such that $g^m$ is conjugate to a periodically geodesic element $h$.
(Moreover, the smallest such integer $m$ is uniformly bounded.)
The rationality of the translation numbers in a word hyperbolic group
is an immediate consequence of this result because the translation
numbers are constant on the conjugacy class and
$\t_X(h^n)=|n|\t_X(h)$ for all elements $h$ and integers $n$.

Before stating our results, we briefly review Garside groups.
The class of Garside groups, first introduced by Dehornoy
and Paris~\cite{DP99}, provides a lattice-theoretic
generalization of braid groups and Artin groups of finite type.
Garside groups are equipped with a special element $\Delta$,
called the Garside element, and a finite set $\D$,
called the set of simple elements.
Elements in a Garside group admit a unique normal form
$$\Delta^r s_1\cdots s_k,
$$
where $r\in\Z$ and $s_1,\ldots,s_k\in\D\setminus\{ 1, \Delta\}$.
There are integer-valued invariants $\inf$ and $\sup$
for elements in a Garside group such that
if $\Delta^r s_1\cdots s_k$ is the normal form of $g$, then
$\inf(g)=r$ and $\sup(g)=r+k$.
For an element $g\in G$, let $[g]$ denote its conjugacy class,
that is, $[g]=\{h^{-1}gh:h\in G\}$.
The conjugacy invariants $\infs$ and $\sups$ are defined as
$\infs(g)=\max\{\inf(h):h\in[g]\}$ and
$\sups(g)=\min\{\sup(h):h\in[g]\}$.
The super summit set of $g$ is defined as
$$[g]^S=\{h\in[g]: \inf(h)=\infs(g),\ \sup(h)=\sups(g)\}.$$
Intuitively, the super summit set is the set
of all elements in the conjugacy class
that have the shortest normal form in that class.
(Note that the number of simple elements $s_i$ in the normal form of $g$
is $\sup(g)-\inf(g)$, which decreases as $\inf(g)$ increases
and as $\sup(g)$ decreases.)

\def\aa{Garside groups are equipped with a special element $\Delta$,
called the Garside element, and a finite set $\D$,
called the set of simple elements.
In the conjugacy class of an element $g$ of a Garside group,
there is a unique finite nonempty subset, called the super summit set.}
We now state the main result of this note.
See the next section for the definition of
$\Vert\cdot\Vert$.

\medskip\noindent
\textbf{Theorem A} (Theorem~\ref{thm:straight})\ \
Let $G$ be a Garside group with Garside element $\Delta$ and the set $\D$ of
simple elements.
Let $\|\Delta\| = N$ and $g\in G$.
There exists a positive integer $n \le N^2$ such that
every element of the super summit set of $g^n$
is periodically geodesic with respect to $\D$.
\medskip

We remark that our approach in Garside groups is different from
Swenson's in word hyperbolic groups~\cite{Swe95}.
Swenson showed the existence of periodically geodesic power
up to conjugacy, which implies the rationality of translation numbers.
For Garside groups, we use the rationality
of translation numbers, which is the main result of~\cite{LL06b},
in order to prove Theorem A,
the existence of periodically geodesic power up to conjugacy.
It seems very difficult to prove Theorem A without using
the rationality of translation numbers.
Therefore, the two results, the rationality of translation numbers
and the existence of periodically geodesic powers,
in word hyperbolic groups and in Garside groups
are established in reverse order.

\medskip
We explain our approach briefly.
It is known that for an element $g$ of a Garside group,
the shortest word length $|g|_\D$ is either
$-\inf(g)$, $\sup(g)$ or $\sup(g)-\inf(g)$.
Because $\sup(g)=-\inf(g^{-1})$, the proof of Theorem A
can be reduced to proving that $g^n$ is conjugate to
an inf-straight element for some $n\le N$.
(An element $h$ is inf-straight if $\inf(h^m)=m\inf(h)$
for all positive integers $m$.)
This is a consequence of the results on the rationality of translation numbers
in~\cite{LL06b} and the inequality established in Theorem~\ref{thm:INF-main}:
$$ \infs(g)\le\INF(g)<\infs(g)+1$$
where $\INF(g)=\lim_{n\to\infty}\inf(g^n)/n$.

\medskip

Another interest of this note is to construct an algorithm
for computing stable super summit sets in Garside groups.
The \emph{stable super summit set} of an element $g$ of a Garside group
is defined as
$$
[g]^\St=\{h\in[g]^S:
\mbox{$h^n\in[g^n]^S$ for all $n\ge 1$} \}.
$$
Namely, it is the set of all conjugates $h$ of $g$
whose every power $h^n$ has the shortest normal form
in the conjugacy class of $g^n$.
It is known that every stable super summit set is nonempty
and satisfies most of the important properties of super summit sets~\cite{LL06a}.
However, as noted in~\cite{LL06a},
we need an affirmative answer to the following question
in order to compute stable super summit sets.
\begin{quote}
Can we make a finite-time algorithm that decides,
given an element $h$ in the conjugacy class of $g$,
whether it belongs to the stable super summit set of $g$?
\end{quote}
In this note, we solve the above problem, and hence obtain the following result.
\par\medskip\noindent
\textbf{Theorem B} (Corollary~\ref{thm:alg_stSSS}).\ \
There is a finite-time algorithm that, given an element
of a Garside group, computes its stable super summit set.

\medskip

\section{Garside groups}
We start with the definition of Garside groups
and known solutions to the conjugacy problem in Garside groups.
See~\cite{Gar69,BKL98,DP99,Deh02,FG03,Pic01b,Geb05}
for details.

\subsection{Garside monoids and groups}

Let $M$ be a monoid.
Let \emph{atoms}  be the elements $a\in M\setminus \{1\}$
such that $a=bc$ implies either $b=1$ or $c=1$.
Let $\Vert a\Vert$ be the supremum of the lengths of all expressions of
$a$ in terms of atoms. The monoid $M$ is said to be \emph{atomic}
if it is generated by its atoms and $\Vert a\Vert<\infty$ for any $a\in M$.
In an atomic monoid $M$, there are partial orders $\le_L$ and $\le_R$:
$a\le_L b$ if $ac=b$ for some $c\in M$;
$a\le_R b$ if $ca=b$ for some $c\in M$.

\begin{definition}
An atomic monoid $M$ is called a \emph{Garside monoid} if
\begin{enumerate}
\item[(i)] $M$ is finitely generated;
\item[(ii)] $M$ is left and right cancellative;
\item[(iii)] $(M,\le_L)$ and $(M,\le_R)$ are lattices;
\item[(iv)] there exists an element $\Delta$, called a
\emph{Garside element}, satisfying the following:\\
(a) for each $a\in M$, $a\le_L\Delta$ if and only if $a\le_R\Delta$;\\
(b) the set $\{a\in M: a \le_L\Delta\}$ generates $M$.
\end{enumerate}
\end{definition}

An element $a$ of $M$ is called a \emph{simple element} if $a\le_L\Delta$.
Let $\D$ denote the set of all simple elements.
Let $\wedgeL$ and $\veeL$ denote the lattice operations
of the poset $(M,\le_L)$.

Garside monoids satisfy Ore's conditions,
and thus embed in their groups of fractions.
A \emph{Garside group} is defined as the group of fractions
of a Garside monoid.
When $M$ is a Garside monoid and $G$ the group of fractions of $M$,
we identify the elements of $M$ and their images in $G$
and call them \emph{positive elements} of $G$.
$M$ is called the \emph{positive monoid} of $G$,
often denoted by $G^+$.

The partial orders $\le_L$ and $\le_R$, and thus the lattice structures
in the positive monoid $G^+$ can be extended
to the Garside group $G$ as follows:
$g\le_L h$ (respectively, $g\le_R h$) for $g,h\in G$
if $gc=h$ (respectively, $cg=h$) for some $c\in G^+$.

For $g\in G$, there are integers $r\le s$ such that
$\Delta^r\le_L g\le_L\Delta^s$.
Hence, the invariants
$\inf(g)=\max\{r\in\Z:\Delta^r\le_L g\}$,
$\sup(g)=\min\{s\in\Z:g\le_L \Delta^s\}$ and
$\len(g)=\sup(g)-\inf(g)$
are well-defined.
It is known that, for $g\in G$, there exists a unique expression
$$
g=\Delta^r s_1\cdots s_k,
$$
such that $s_1,\ldots,s_k\in \D\setminus\{1,\Delta\}$,
$(s_is_{i+1}\cdots s_k)\wedgeL \Delta=s_i$ for $i=1,\ldots,k$,
$\inf(g)=r$, and\/ $\sup(g)=r+k$.
Such an expression is called the \emph{normal form} of $g$.

For $g\in G$, we denote its conjugacy class $\{ h^{-1}gh : h\in G\}$
by $[g]$.
Define $\infs(g)=\max\{\inf(h):h\in [g]\}$ and
$\sups(g)=\min\{\sup(h):h\in [g]\}$.
The \emph{super summit set} $[g]^S$
and the \emph{stable super summit set} $[g]^\St$
are subsets of the conjugacy class of $g$ defined as follows:
\begin{eqnarray*}
[g]^S&=&\{h\in [g]:\inf(h)=\infs(g) \ \mbox{ and } \sup(h)=\sups(g)\};\\{}
[g]^\St&=&\{h\in [g]^S:h^k\in[g^k]^S \ \mbox{ for all positive integers $k$}\}.
\end{eqnarray*}
It is well known that $[g]^S$ is finite and nonempty.
Since $[g]^\St$ is a subset of $[g]^S$,
$[g]^\St$ is also a finite set.
It is proved in \cite{BGG06a,LL06a}
that $[g]^\St$ is nonempty.

In the rest of the paper, if it is not specified,
$G$ is assumed to be a Garside group, whose positive monoid is $G^{+}$,
with Garside element $\Delta$ and the set $\D$ of simple elements,
where $\| \Delta\|$ is simply written as $N$.

\subsection{Some results on $\inf$ and $\sup$}
Here, we collect some results on $\inf$ and $\sup$.
See \cite{LL06a} for Lemma~\ref{thm:Lmax} and Theorem~\ref{thm:inequality},
and \cite{LL06b} for the others.

\medskip

For a positive element $a\in G^+$, define $\Lmax(a)=\Delta\wedge_L a$.

\begin{lemma}\label{thm:Lmax}
Let $g\in G$ and $a\in G^+$.
If\/ $\inf(ga)>\inf(g)$, then $\inf(g\Lmax(a))>\inf(g)$.
\end{lemma}

\begin{theorem}\label{thm:inequality}
Let $g\in G$.
For all $n\ge 1$, the following hold.
\begin{enumerate}
\item[(i)]
$n\infs(g)\le \infs(g^n)\le n\infs(g)+n-1$.
\item[(ii)]
$n\sups(g)-(n-1)\le\sups(g^n)\le n\sups(g)$.
\end{enumerate}
\end{theorem}

It is proved in~\cite[Lemmas 3.2 and 3.3]{LL06b} that
the following limits are well-defined for all elements $g\in G$:
$$
\INF(g)=\lim_{n\to\infty}\frac{\inf(g^n)}n;\quad
\SUP(g)=\lim_{n\to\infty}\frac{\sup(g^n)}n;\quad
\LEN(g)=\lim_{n\to\infty}\frac{\len(g^n)}n.
$$

\begin{lemma}\label{thm:INF-basic}
Let $g,h\in G$.
\begin{enumerate}
\item[(i)] $\INF(h^{-1}gh)=\INF(g)$ and $\SUP(h^{-1}gh)=\SUP(g)$.
\item[(ii)] For all $n\ge 1$, $\INF(g^n)=n\INF(g)$ and $\SUP(g^n)=n\SUP(g)$.
\item[(iii)] $\infs(g) \le \INF(g) \le \infs(g)+1$ and
$\sups(g)-1 \le \SUP(g) \le \sups(g)$.
\item[(iv)] Both $\INF(g)$ and $\SUP(g)$ are rational numbers of
the form $p/q$ for some integers $p$ and $q$ with $1\le q\le N$.
\end{enumerate}
\end{lemma}

For an element $g\in G$,
\begin{enumerate}
\item[(i)] $g$ is said to be \emph{inf-straight} if $\inf(g)=\INF(g)$;
\item[(ii)]
$g$ is said to be \emph{sup-straight} if $\sup(g)=\SUP(g)$.
\end{enumerate}

\begin{lemma}\label{thm:inf-straight-equiv}
For every $g\in G$, the following conditions are equivalent.
\begin{enumerate}
\item[(i)] $g$ is inf-straight.
\item[(ii)] $\inf(g^N)=N\inf(g)$.
\item[(iii)] $\inf(g^n)=n \inf(g)$ for all $n\ge 1$.
\end{enumerate}
\end{lemma}

\begin{lemma}\label{thm:sup-straight-equiv}
For every $g\in G$, the following conditions are equivalent.
\begin{enumerate}
\item[(i)] $g$ is sup-straight.
\item[(ii)] $\sup(g^N)=N\sup(g)$.
\item[(iii)] $\sup(g^n)=n \sup(g)$ for all $n \ge 1$.
\end{enumerate}
\end{lemma}

\begin{lemma}\label{thm:inf-straight-conj}
For every $g\in G$, the following conditions are equivalent.
\begin{enumerate}
\item[(i)]
$g$ is conjugate to an inf-straight element.
\item[(ii)]
$\infs(g)=\INF(g)$.
\item[(iii)]
$\infs(g^N)=N\infs(g)$.
\item[(iv)]
$\infs(g^k)=k\infs(g)$ for all $k\ge 1$.
\item[(v)]
For all $h\in[g]^S$, $h$ is inf-straight.
\end{enumerate}
\end{lemma}

\begin{lemma}\label{thm:sup-straight-conj}
For every $g\in G$, the following conditions are equivalent.
\begin{enumerate}
\item[(i)]
$g$ is conjugate to a sup-straight element.
\item[(ii)]
$\sups(g)=\SUP(g)$.
\item[(iii)]
$\sups(g^N)=N\sups(g)$.
\item[(iv)]
$\sups(g^k)=k\sups(g)$ for all $k\ge 1$.
\item[(v)]
For all $h\in[g]^S$, $h$ is sup-straight.
\end{enumerate}
\end{lemma}

\section{Asymptotic limit of $\inf$}

For a real number $x$, let $\lfloor x\rfloor$ denote
the largest integer less than or equal to $x$.
From Lemma~\ref{thm:INF-basic}~(iii),
$\infs(g)\le \INF(g)\le\infs(g)+1$
for all elements $g$ of a Garside group.
The goal of this section is to show that
$\INF(g)=\infs(g)+1$ cannot happen,
hence $\infs(g)=\lfloor\INF(g)\rfloor$.

We first recall Schur's theorem.
Schur~\cite{Sch16} proved that for every positive integer $M$,
there exists a positive integer $L$ such that for every partition
of the set $\{1,2,\ldots,L\}$ into $M$ subsets,
one of the subsets contains two numbers $n$ and $m$ together
with their sum $n+m$.
The smallest such integer $L$ is called the {\em Schur number}.
As a corollary, we have the following lemma.

\begin{lemma}\label{lem:Shur}
Let $\N$ be the set of all positive integers.
Let $T_1,\ldots, T_M$ be a finite collection of
subsets of\/ $\N$ such that $\N=T_1\cup\cdots\cup T_M$.
Then there is at least one set $T_k$ that contains
$n$, $m$ and $n+m$
for some $n, m\in\N$.
\end{lemma}

\begin{theorem}\label{thm:INF-main}
For every element $g$ of a Garside group $G$,
$\infs(g)\le\INF(g)<\infs(g)+1$.
\end{theorem}

\begin{proof}
Because $\infs(g)\le\INF(g)\le\infs(g)+1$ by Lemma~\ref{thm:INF-basic}~(iii),
it suffices to show that $\INF(g)\ne\infs(g)+1$.
On the contrary, assume
\begin{equation}\label{assumption}
\INF(g)=\infs(g)+1.
\end{equation}
Since both $\INF$ and $\infs$ are conjugacy invariants, we may
assume that $g$ belongs to its stable super summit set, and
hence
$$
\inf(g^n)=\infs(g^n)\quad\mbox{and}\quad
\sup(g^n)=\sups(g^n)
 \ \ \mbox{for all}\ n\ge 1.$$

\begin{claim}\label{claim1}
For all $n,m\ge 1$, the following hold:
\begin{enumerate}
\item[(i)]   $\inf(g^n)=n\inf(g)+n-1$;
\item[(ii)]  $\inf(g^{n+m})=\inf(g^n)+\inf(g^m)+1$;
\item[(iii)] $\len(g^n)\ge 1$.
\end{enumerate}
\end{claim}

\begin{proof}[Proof of Claim~\ref{claim1}]
Recall that $\inf(g^n)=\infs(g^n)$ and $\sup(g^n)=\sups(g^n)$
for all $n\ge 1$, because we have
assumed that $g$ belongs to its stable
super summit set.

(i)\ \
By Theorem~\ref{thm:inequality}~(i), $\infs(g^n)\le n\infs(g)+n-1$ for all $n\ge 1$.
Assume that $\infs(g^n)\le n\infs(g)+n-2$ for some $n\ge 1$.
By Lemma~\ref{thm:INF-basic}~(iii),
$\INF(g^n)\le \infs(g^n)+1\le n\infs(g)+n-1$, and thus
$$
\INF(g)
=\frac{\INF(g^n)}n
\le\infs(g)+1-\frac1n
<\infs(g)+1.
$$
This contradicts the assumption (\ref{assumption})
which states $\INF(g)=\infs(g)+1$.

\smallskip(ii)\ \
By (i),
$\begin{array}[t]{l}
\infs(g^n)+\infs(g^m)+1 =(n\infs(g)+n-1)+(m\infs(g)+m-1)+1\\
\mbox{}\qquad =(n+m)\infs(g)+n+m-1=\infs(g^{n+m}).
\end{array}$

\smallskip(iii)\ \
If $\lens(g^n)=0$ for some $n\ge 1$,
then $g^n$ is conjugate to $\Delta^k$ for some integer $k$.
Therefore $\infs(g^{2n})=2k=2\infs(g^n)$ is an even integer.
On the other hand, $\infs(g^{2n})=2n\infs(g)+2n-1$ by (i),
hence $\infs(g^{2n})$ is an odd integer.
It is a contradiction.
\end{proof}

For $n\ge 1$, let $r_n=\inf(g^n)$ and
$s_n=\Lmax(g^n\Delta^{-r_n})$.
Then,
$$ g^n=s_na_n\Delta^{r_n}, $$
for some $a_n\in G^+\backslash\Delta G^{+}$.
By Claim~\ref{claim1}, for all $n,m\ge 1$,
$s_n\ne 1$ and
\begin{equation}
r_{n+m}=r_n+r_m+1.\label{eqn:inf_n+m}
\end{equation}

\begin{claim}\label{claim2}
If\/ $s_n=s_{n+m}$ for some $n, m\ge 1$, then $\inf(a_n\Delta^{r_n}s_m)=r_n+1$.
\end{claim}

\begin{proof}[Proof of Claim~\ref{claim2}]
Note that
$$
(s_na_n\Delta^{r_n}) (s_ma_m\Delta^{r_m})
=g^ng^m
=g^{n+m}
=s_{n+m}a_{n+m}\Delta^{r_{n+m}}.
$$
Since $s_n=s_{n+m}$ and $r_{n+m}=r_n+r_m+1$,
we obtain $a_n\Delta^{r_n}s_ma_m=a_{n+m}\Delta^{r_n+1}$.
Therefore
$$
\inf(a_n\Delta^{r_n}s_ma_m)
= r_n+1 >\inf(a_n\Delta^{r_n}).
$$
Since $s_m=\Lmax(s_ma_m)$,
$\inf(a_n\Delta^{r_n}s_m)>\inf(a_n\Delta^{r_n})=r_n$
by Lemma~\ref{thm:Lmax}.
Therefore,
$$\inf(a_n\Delta^{r_n}s_m)\ge r_n+1.$$
On the other hand,
$\inf(a_n\Delta^{r_n}s_m) \le\inf(a_n\Delta^{r_n})+1=r_n+1$.
As a consequence, we obtain $\inf(a_n\Delta^{r_n}s_m)=r_n+1$.
\end{proof}

For each simple element $s\in\D$, let
$T_s=\{n\in\N: s_n=s\}$.
Note that $\N=\cup_{s\in\D} T_s$
and $\D$ is a finite set.
By Lemma~\ref{lem:Shur}, there exists a subset $T_s$ that contains
$n$, $m$ and $n+m$ for some $n,m\ge 1$.
($n$ and $m$ may not be distinct.)
By definition of $T_s$,
$$
s_n=s_m=s_{n+m}.
$$
Applying Claim~\ref{claim2} to $(n,n+m)$ and $(m,n+m)$,
we obtain
$$
\inf(a_n\Delta^{r_n}s_m)=r_n+1\quad\mbox{and}\quad
\inf(a_m\Delta^{r_m}s_n)=r_m+1,
$$
from which
\begin{eqnarray*}
r_{n+m}&=&\infs(g^{n+m}) \ge \inf(s_n^{-1} g^{n+m}s_n)
=\inf(s_n^{-1}g^ng^ms_n) \\
 &=& \inf(s_n^{-1}(s_na_n\Delta^{r_n}) (s_ma_m\Delta^{r_m}) s_n)
  = \inf((a_n\Delta^{r_n}s_m)(a_m\Delta^{r_m}s_n)) \\
& \ge & \inf(a_n\Delta^{r_n}s_m)+\inf(a_m\Delta^{r_m}s_n)
= (r_n+1)+(r_m+1) \\
& = & r_n+r_m+2.
\end{eqnarray*}
\def\temp{
Then $\infs(g^{n+m})\ge r_n+r_m+2$ because
\begin{eqnarray*}
\lefteqn{\inf(s_n^{-1} g^{n+m}s_n)
=\inf(s_n^{-1}g^ng^ms_n)
=\inf(s_n^{-1}(s_na_n\Delta^{r_n}) (s_ma_m\Delta^{r_m}) s_n)}\\
&=& \inf((a_n\Delta^{r_n}s_m)(a_m\Delta^{r_m}s_n))
\ge \inf(a_n\Delta^{r_n}s_m)+\inf(a_m\Delta^{r_m}s_n)\\
&=& (r_n+1)+(r_m+1).
\end{eqnarray*}
}
This contradicts (\ref{eqn:inf_n+m}).
\end{proof}

\begin{corollary}\label{cor:TBound_tinf}
Let $G$ be a Garside group with Garside element $\Delta$,
and let $N=\Vert\Delta\Vert$.
For every element $g\in G$, we have
\begin{itemize}
\item[(i)] $\infs(g)\le\INF(g)\le\infs(g)+1-1/N$;
\item[(ii)] $\sups(g)-1+1/N\le\SUP(g)\le\sups(g)$;
\item[(iii)] $\lens(g)-2+2/N\le\LEN(g)\le\lens(g)$.
\end{itemize}
\end{corollary}

\begin{proof}
(i)\ \
Since $\INF(g)$ is of the form $p/q$ for some $p,q\in\Z$ with
$1\le q\le N$ by Lemma~\ref{thm:INF-basic}~(iv),
and $\infs(g)\le\INF(g)<\infs(g)+1$
by Theorem~\ref{thm:INF-main},
we have $\infs(g)\le\INF(g)\le\infs(g)+1-1/N$ as desired.

(ii)\ \
It follows from (i) because $\sups(g)=-\infs(g^{-1})$.

(iii)\ \
It follows from (i) and (ii).
\end{proof}

The following example shows that
the upper bound $\infs(g)+1-1/N$
of $\INF(g)$ in Corollary~\ref{cor:TBound_tinf}~(i)
is optimal.

\begin{example}\label{eg:1}
For an integer $N\ge 2$, let
$$G=\langle x,y\mid x^N=y^N\rangle.$$
It is a Garside group with Garside element
$\Delta=x^N=y^N$~\cite[Example 4]{DP99},
and $\Vert\Delta\Vert=N$.
Let $g=x^{N-1}$. Then $\inf(g^k)=\lfloor k(N-1)/N\rfloor$
for all $k\ge 1$,
hence $\INF(g)=1-1/N$. Since $\infs(g)=0$,
$\INF(g)=\infs(g)+1-1/N$.
\long\def\temp{
Let $g_1=x^{N-1}$. Then $\inf(g_1^k)=\lfloor k(N-1)/N\rfloor$
for all $k\ge 1$,
hence $\INF(g_1)=1-1/N$. Since $\infs(g_1)=0$,
$\INF(g_1)=\infs(g_1)+1-1/N$.

Similarly, if we take $g_2=x$, then $\sups(g_2^k)=\lceil k/N\rceil$,
hence $\SUP(g_2)=1/N$.
Since $\sups(g_2)=1$, $\SUP(g_2)=\sups(g_2)-1+1/N$.
}
\end{example}

\def\temp{
\begin{example} Delete??? \\
Let $H=\langle x,y\mid x^p=y^q\rangle$, $p=q+1\ge 2$.
Let $G=H\times H$ and $g=(x^{-1},y)$.
Then $\inf(g^k)=\lfloor -k/p\rfloor$
and $\sup(g^k)=\lceil k/q\rceil$.
Therefore $\len(g^k)=\lceil k/q\rceil-\lfloor -k/p\rfloor$, hence
$$\LEN(g)=\lim_{n\to\infty} \len(g^k)/k
 =1/q-(-1/p)=(p+q)/pq$$
Note that $\lens(g)=1$. Hence
$\LEN(g)-(\lens(g)-2)=(p+q)/pq=(4N-1)/N(N-2)$ is asymptotically $4/N$.
\end{example}
}

From Corollary~\ref{cor:TBound_tinf} and
Lemmas~\ref{thm:inf-straight-conj} and~\ref{thm:sup-straight-conj},
we have the following.
\begin{corollary}\label{thm:inf_sup_str}
Let $g$ be an element of a Garside group.
\begin{itemize}
\item[(i)] $g$ is conjugate to an inf-straight element
if and only if $\INF(g)$ is an integer.
\item[(ii)] $g$ is conjugate to a sup-straight element
if and only if $\SUP(g)$ is an integer.
\end{itemize}
\end{corollary}

By Lemma~\ref{thm:inf-straight-conj} and Corollary~\ref{thm:inf_sup_str},
$\INF(g)$ is an integer
if and only if $\infs(g^n)=n\infs(g)$ for all $n\ge 1$,
and the same is true for $\SUP(g)$ and $\sups(g)$.
However, the following example shows that
we cannot expect such a property for $\LEN(g)$.

\begin{example}
Consider the group
$$
G=H\times H,\quad
\mbox{where $H=\langle x,y\mid x^{2p}=y^{2p}\rangle$ for $p\ge 1$.}
$$
As in Example~\ref{eg:1},
$H$ is a Garside group with Garside element $\Delta_H =x^{2p}=y^{2p}$.
Because $G$ is a cartesian product of $H$,
it is a Garside group with Garside element
$\Delta = (\Delta_H, \Delta_H)$ by~\cite[Theorem 4.1]{Lee06}.
Let $g=(x^{-p},y^{p})$.
Then, for all $n\ge 1$, $g^n = (x^{-np}, y^{np})$ and
it is not difficult to show that
$$
\inf(g^n)=\infs(g^n)=\lfloor -n/2\rfloor \quad\mbox{and}\quad
\sup(g^n)=\sups(g^n)=\lceil n/2\rceil.
$$
Therefore $\INF(g)=-1/2$ and $\SUP(g)=1/2$
and it follows that $\LEN(g)=\SUP(g)-\INF(g)=1$.
In particular $\LEN(g)$ is an integer.
However, it is not true that $\lens(g^n)=n\lens(g)$ for all $n\ge 1$.
(For example, $\lens(g^2)=\lens(g)=2$.)
\end{example}

\section{Periodically geodesic elements}

\begin{definition}
Let $G$ be a group and $X$ be a finite set of semigroup generators for $G$.
An element $g\in G$ is said to be \emph{periodically geodesic}
with respect to $X$
if $|g^n|_X=|n|\cdot |g|_X$ for all $n\in\Z$.
\end{definition}

Note that if $X$ is closed under inversion, then $|g|_X=|g^{-1}|_X$,
hence $g$ is periodically geodesic
if $|g^n|_X=n |g|_X$ for all $n\ge 1$.

\begin{lemma}\label{lem:straight}
Let $G$ be a Garside group with the set $\D$ of simple elements, and
let $g\in G$.
Every element $h\in [g]^S$ is periodically geodesic with respect to $\D$
if one of the following conditions holds.
\begin{itemize}
\item[(i)] $\infs(g)\ge 0$ and $\SUP(g)$ is an integer.
\item[(ii)] $\sups(g)\le 0$ and $\INF(g)$ is an integer.
\item[(iii)] Both $\INF(g)$ and $\SUP(g)$ are integers.
\end{itemize}
\end{lemma}

\begin{proof}
It is well known (see~\cite{Cha95} for example) that
for any element $g\in G$
$$|g|_{\D}=\left\{\begin{array}{rl}
\sup(g) & \mbox{if } \inf(g)\ge 0; \\
-\inf(g) & \mbox{if } \sup(g)\le 0; \\
\len(g) & \mbox{otherwise.}
\end{array}
\right.
$$

(i)\ \
Since $\infs(g)\ge 0$, $\inf(h^n) \ge n\inf(h)=n\infs(g)\ge 0$
for all $n\ge 1$,
hence
$$|h^n|_\D=\sup(h^n) \quad\mbox{for all}\quad n\ge 1.$$
Since $\SUP(g)$ is an integer,
$g$ is conjugate to a sup-straight element
by Corollary~\ref{thm:inf_sup_str}~(ii).
Since $h\in [g]^S$,
$h$ is sup-straight by Lemma~\ref{thm:sup-straight-conj}.
Thus, by Lemma~\ref{thm:sup-straight-equiv}
$$\sup(h^n)=n\sup(h) \quad\mbox{for all}\quad n\ge 1.$$

Consequently, $|h^n|_\D=\sup(h^n)=n\sup(h)=n|h|_\D$ for all $n\ge 1$.

\smallskip(ii)\ \
It can be proved similarly to (i).

\smallskip(iii)\ \
If $\infs(g)\ge 0$ or $\sups(g)\le 0$,
then $h$ is periodically geodesic by (i) and (ii).
So, we may assume $\infs(g)<0<\sups(g)$.
Then $\inf(h^n)\le\infs(g^n)<0<\sups(g^n)\le\sup(h^n)$ for all $n\ge 1$
by Theorem~\ref{thm:inequality}, hence
$$| h^n |_\D =\len(h^n) \quad\mbox{for all}\quad n\ge 1.$$
Since both $\INF(g)$ and $\SUP(g)$ are integers,
$g$ is conjugate to an element which is inf-straight and
to an element which is sup-straight, by Corollary~\ref{thm:inf_sup_str}.
Since $h\in [g]^S$,
$h$ is both inf- and sup-straight by
Lemmas~\ref{thm:inf-straight-conj} and ~\ref{thm:sup-straight-conj}.
Then, by Lemmas~\ref{thm:inf-straight-equiv}
and~\ref{thm:sup-straight-equiv}
$$\len(h^n)=n\len(h) \quad\mbox{for all}\quad n\ge 1.$$

Consequently, $| h^n |_\D =\len(h^n)=n\len(h)=n |h|_\D$
for all $n\ge 1$.
\end{proof}

We now establish the main theorem of this note.

\begin{theorem}\label{thm:straight}
Let $G$ be a Garside group with Garside element $\Delta$ and the set $\D$ of
simple elements.
Let $\|\Delta\| = N$ and $g\in G$.
There exists a positive integer $n \le N^2$ such that
every element of the super summit set of $g^n$
is periodically geodesic with respect to $\D$.
\end{theorem}

\begin{proof}
Let $\INF(g)=p_1/q_1$ and $\SUP(g)=p_2/q_2$,
where $p_1, p_2, q_1, q_2\in\Z$ and $1\le q_1, q_2\le N$.
Let $n$ be the least common multiple of $q_1$ and $q_2$.
Then, $1\le n \le N^2$.
Since both
$\INF(g^n)=np_1/q_1$ and $\SUP(g^n)=np_2/q_2$ are integers,
every element of the super summit set of $g^n$
is periodically geodesic with respect to $\D$ by Lemma~\ref{lem:straight}.
\end{proof}

We remark that we can construct a finite-time algorithm that, given $g\in G$,
computes the power $n$ in the above theorem by Lemma~\ref{lem:straight} and
by using the algorithms
for $\INF(\cdot )$ described in~\S\ref{sec:INF}.
(Applying those algorithms for $\INF(\cdot )$,
we can do the same task for $\SUP(\cdot )$.)

\section{Computation of $\INF(\cdot)$}\label{sec:INF}

In~\cite[Theorem 3.9 (iii)]{LL06b}, there is a finite-time algorithm
for computing $\INF(g)$ given an element $g$ of a Garside group.
This algorithm exploits two facts:
\begin{itemize}
\item[(i)] $\INF(g) = p/q$ for some integers $p$, $q$ with $1\le q\le N$.
\item[(ii)] $\infs(g^{n})/n\le \INF(g)\le (\infs(g^{n})+1)/n$
for all $n\ge 1$.
\end{itemize}
Therefore, for any $n\ge N^2$, $\INF(g)$ is the unique rational number
of the form $p/q$ in the closed interval
$[ \infs(g^{n})/n, (\infs(g^{n})+1)/n ]$,
where $p, q\in\Z$ with $1\le q\le N$.
Applying Theorem~\ref{thm:INF-main},
this section shows different methods
of computing $\INF(\cdot)$.

\begin{theorem}\label{thm:INF-computation}
Let $G$ be a Garside group with Garside element $\Delta$,
and let $N=\|\Delta\|$.
For every $g\in G$, the following hold.
\begin{enumerate}
\item[(i)] $\INF(g)=\max\{\infs(g^k)/k : k=1,\ldots, N\}$.
\item[(ii)]
Let $1\le q\le N$.
$\INF(g)=p/q$ for some $p\in\Z$ if and only if\/
$\infs(g^{qN})= N\infs(g^q)$.
Furthermore, $p=\infs(g^q)$.
\end{enumerate}
\end{theorem}

\begin{proof}
(i)\ \
Since $\INF(g)=\INF(g^k)/k\ge\infs(g^k)/k$ for all $k\ge 1$,
it suffices to show that there exists $k\in\{ 1,\ldots, N\}$ such that
$\INF(g)=\infs(g^k)/k$.
By Lemma~\ref{thm:INF-basic}~(iv),
$\INF(g)=p/q$ for some integers $p, q$ with $1\le q\le N$.
By Theorem~\ref{thm:INF-main}, $\infs(g^q)=\lfloor \INF(g^q)\rfloor
= \lfloor q\INF(g)\rfloor=p$, thus
$$\infs(g^q)/q=p/q=\INF(g).$$

\smallskip (ii)\ \
Suppose that $\infs(g^{qN})= N\infs(g^q)$.
Then, $\infs(g^q)= \INF(g^q)= q\INF(g)$
by Lemma~\ref{thm:inf-straight-conj}.

Conversely, suppose that $\INF(g)=p/q$ for some integer $p$.
Then $\INF(g^q)=q\INF(g)=p$ is an integer, and it follows that
$g^q$ is conjugate to an inf-straight element
by Corollary~\ref{thm:inf_sup_str}~(i).
This means that $\infs(g^{qN})= N\infs(g^q)$
by Lemma~\ref{thm:inf-straight-conj}.
\end{proof}

From Theorem~\ref{thm:INF-computation},
we obtain two ways to compute $\INF(g)$.
\begin{itemize}
\item[(i)] Compute $\infs(g^k)$ for $k=1,\ldots,N$,
and then compute the maximum of $\{\infs(g^k)/k : k=1,\ldots,N \}$.
Then, $\INF(g)$ is this maximum value.
\item[(ii)]
Compute $\infs(g^k)$ and $\infs(g^{kN})$ for $k=1,\ldots,N$,
and find $1\le q\le N$ with $\infs(g^{qN}) = N\infs(g^q)$.
Then, $\INF(g)=\infs(g^q)/q$.
\end{itemize}

\def\temp{
Therefore, we have the following corollary.

\begin{corollary}
There is a finite-time algorithm that,
given an element $g$ of a Garside group,
computes $\INF(g)$.
\end{corollary}
}

\section{Computation of stable super summit sets}

Recall the definition of the stable super summit set $[g]^\St$ of $g\in G$:
$$[g]^\St=\{h\in[g]^S:
h^n\in[g^n]^S \ \mbox{for all}\  n\ge 1 \}.$$
As noted in \S1, in order to compute the stable super summit set
$[g]^\St$,
we need a finite-time algorithm that decides,
given an element $h$ in the conjugacy class of $g$,
whether $h$ belongs to $[g]^\St$.
The existence of such an algorithm is obvious for
the super summit set and the ultra summit set.
However, for the stable super summit set, a naive algorithm will test
whether $h^n\in[g^n]^S$ for all $n\ge 1$.
This kind of algorithm does not halt in finite time.
We resolve this problem by using Theorem~\ref{thm:INF-main}.

\begin{theorem}\label{thm:UBforStableSSS}
Let $G$ be a Garside group with Garside element $\Delta$.
Let $N=\|\Delta\|$ and $g\in G$.
\begin{enumerate}
\item[(i)]
If\/ $h\in [g]$ satisfies
$\inf(h^n)=\infs(g^n)$ for all $1\le n\le N$, then
$\inf(h^n)=\infs(g^n)$ for all $n \ge 1$.

\item[(ii)]
If\/ $h\in [g]$ satisfies
$\sup(h^n)=\sups(g^n)$ for all $1\le n\le N$, then
$\sup(h^n)=\sups(g^n)$ for all $n \ge 1$.

\item[(iii)]
$[g]^\St=\{h\in[g]^S:
h^n \in[g^n]^S \ \mbox{for all}\  1\le n\le N\}$.
\end{enumerate}
\end{theorem}

\begin{proof}
We prove only (i), because (ii) can be proved similarly
and (iii) is a consequence of (i) and (ii).
Recall that $\INF(g)=p/q$ for some integers $p,q$ with $1\le q\le N$.
From Theorem~\ref{thm:INF-main},
$$\infs(g^n)=\lfloor\INF(g^n)\rfloor = \lfloor n \INF(g)\rfloor = \lfloor np/q \rfloor
\quad\mbox{for all}\  n \ge 1.$$
Choose any integer $n\ge 1$.
Then, there exist non-negative integers $k$ and $r$ with $r<q$
such that $n=kq+r$.
Since $\inf(h^i)=\infs(g^i)$ for all $1\le i \le N$ by the hypothesis,
\begin{eqnarray*}
\inf(h^n)
&\ge& k\inf(h^q)+\inf(h^r)
=k\infs(g^q)+\infs(g^r)
=kp+\lfloor rp/q\rfloor\\
&=&\lfloor (kq+r)p/q\rfloor
= \lfloor n \INF(g)\rfloor
= \infs(g^n).
\end{eqnarray*}
By the definition of $\infs$, $\inf(h^n)\le\infs(g^n)$.
Consequently, $\inf(h^n)=\infs(g^n)$.
Since $n$ was arbitrarily chosen, $\inf(h^n)=\infs(g^n)$ for all $n\ge 1$.
\end{proof}

From Theorem~\ref{thm:UBforStableSSS}~(iii),
it suffices to check if $h^n \in [g^n]^S$ only for all $1\le n \le N$
in order to decide $h\in [g]^\St$.

\begin{corollary}\label{cor:StSSSdecision}
There is a finite-time algorithm that, given elements $g$ and $h$
of a Garside group, decides whether $h\in[g]^\St$ or not.
\end{corollary}

Using the algorithm in the above corollary, together with
the results in \cite{LL06a}, we can make a finite-time
algorithm for computing the stable super summit set:

\begin{itemize}
\item[(i)]
{\em Compute an element $h_0$ in $[g]^\St$.} \\
By Corollary 3.12 in~\cite{LL06a},
there is a finite-time algorithm that, given $g\in G$ and $n\ge 1$,
computes an element $h\in [g]$ with the property that
$h^k\in[g^k]^S$ for all $k\in\{ 1,\ldots,n\}$.
By Theorem~\ref{thm:UBforStableSSS}, such an element $h$ belongs
to $[g]^\St$ if we take $n\ge N$.

\item[(ii)]
{\em Compute all the elements in $[g]^\St$
from $h_0\in[g]^\St$ obtained in the above step.} \\
By Corollary 4.5~(iii) in~\cite{LL06a}, for any $h\in[g]^\St$,
there exists a finite sequence
$$h_0\to h_1\to\cdots\to h_m=h$$
such that for each $i=1,\ldots,m$, $h_i\in[g]^\St$ and $h_i=s_i^{-1}h_{i-1}s_i$
for some simple element $s_i$.
Therefore, using Corollary~\ref{cor:StSSSdecision},
the stable super summit set can be computed in the way to compute
super summit sets:
Define $V_1=\{h_0\}$, and then recursively compute
$V_i=\{ s^{-1}hs: s\in\D, h\in V_{i-1}\}\cap[g]^\St$;
Then $V_1\subset V_2\subset\cdots\subset [g]^\St$;
Since $[g]^\St$ is a finite set, there exists $k\ge 1$ such that
$V_k=V_{k+1}$;
Then $V_k=[g]^\St$.
\end{itemize}

\begin{corollary}\label{thm:alg_stSSS}
There is a finite-time algorithm that, given an element $g$
of a Garside group, computes $[g]^\St$.
\end{corollary}

\end{document}